\documentclass[12pt]{amsart}
\usepackage{amsfonts}
\usepackage{amsmath}
\usepackage{amscd}
\usepackage{amssymb}
\usepackage[dvips]{graphicx}
\usepackage[dvips]{epsfig}
\usepackage{color}
\def\?{\char'76}
\def\!{\char'74}

\def\s{\mathbb{S}^1}

\def\Im{\text{\rm Im}}

\def\L{\mathcal{L}}

\def\q{\overline{q}}

\def\G#1#2{\mathcal{G}_{_{\q}}^{^{#1}}(#2)}

\def\Resizq#1#2{\mathfrak{Upp}_{_{\q}}^{^{#1}}(#2)}

\newcounter{numero}

\newcounter{letra}

\newcounter{romnumero}

\newcounter{bibnumero}

\theoremstyle{definition}

\theoremstyle{remark}

\begin{document}
\pagestyle{myheadings} \markboth{G. Padilla}{Intersection
Cohomology of $\s$-Actions on Pseudomanifolds}

\title{Intersection Cohomology of $\s$-Actions on Pseudomanifolds}
\author{G. Padilla}
\address{Universidad Central de Venezuela-
Escuela de Matem\'atica. Caracas 1010.}
\email{gabrielp@euler.ciens.ucv.ve} \dedicatory{ To my wife with
love.}
\date{May 9/2002}
\keywords{Intersection Cohomology, Stratified Pseudomanifolds}
\subjclass{35S35; 55N33}
\begin{abstract}
    For any smooth free action of the unit circle $\s$
    in a manifold $M$; the Gysin sequence of $M$ is a long exact sequence
    relating the DeRham cohomologies of $M$ and its orbit space $M/\s$.
    If the action is not free then $M/\s$ is not a manifold but a stratified
    pseudomanifold and there  is a Gysin
    sequence relating the DeRham cohomology of $M$ with the
    intersection cohomology of $M/\s$. In this work we extend the
    above statements for any stratified pseudomanifold
    $X$ of length 1, whenever the action of $\s$ preserves
    the local structure. We give a Gysin sequence
    relating the intersection cohomologies of $X$ and $X/\s$ with a third
    term $\mathcal{G}$, the Gysin term; whose  cohomology depends
    on basic cohomological data of two flavors: global data concerns the
    Euler class induced by the action, local data relates the Gysin term
    and the cohomology of the fixed strata with values on a locally trivial
    presheaf.
\end{abstract}
\maketitle

\section*{Foreword}

{\it For an entire version of this article the reader should go better to
\tt Indag. Math. Vol. 15 (3), 383-412.}\vskip5mm

A pseudomanifold is a topological space $X$ with two features.
First, there is a closed $\Sigma\subset X$ called the singular
part, which is the disjoint union of smooth manifolds. The
$X-\Sigma$ is a dense smooth manifold. We call strata the
connected components of $\Sigma$ and $X-\Sigma$; they constitute a
locally finite partition of $X$. The second feature is the local
conical behavior of $X$, the model being a product $U\times c(L)$
of a smooth manifold $U$ with the open cone of a compact smooth
manifold $L$ called the link of $U$. A careful reader will notice
that stratified pseudomanifolds with arbitrary length have a
richer and more complicated topological structure; in this article
we deal with stratified pseudomanifolds of length $\leq1$, which
we call just {\it pseudomanifolds}.\newline

Between the various ways for defining the intersection
(co)homology; the reader can see \cite{gm1},\cite{borel2} for a
definition in $pl$-stratified pseudomanifolds;
\cite{borel},\cite{gm2},\cite{pervsheaves} for a definition with
sheaves; \cite{nagase} for an approach with $\L^2$-cohomology;
\cite{brylinsky} for an exposition in Thom-Mather spaces.\newline

In this article, we use the DeRham-like definition exposed in
\cite{illinois} where the  reader will find a beautiful proof of
the DeRham theorem for stratified spaces. We work with
differential forms in $X-\Sigma$ and measure their behavior when
approaching to $\Sigma$, trough an auxiliary construction called
an {\it unfolding} of $X$. Although $X$ may have many different
unfoldings, its intersection cohomology does not depend on any
particular choice. This point of view is the dual of the
intersection  homology defined by King \cite{king}, who works with
a broader family of perversities. When $\s$ acts on $X$ preserving
the local structure then the orbit space $X/\s$ is again a
pseudomanifold with an unfolding. \newline

The well known Gysin sequence of a smooth manifold $M$ with a
principal action of $\s$ is the long exact sequence
\[
    \cdots\rightarrow H^{^i}(M)\overset{\oint}{\rightarrow}
    H^{^{i-1}}(M/\s)\overset{\varepsilon}\rightarrow
    H^{^{i+1}}(M/\s)\overset{\pi^{*}}\rightarrow
    H^{^{i+1}}(M)\rightarrow\cdots
\]
where $\pi^{*}$ is induced by the orbit map
$\pi:M\rightarrow M/\s$, which is a smooth $\s$-principal bundle.
The map $\oint$ is induced by the
integration along the fibers and the connecting homomorphism $\varepsilon$ is
the multiplication by the Euler class $\varepsilon\in
H^2(M/\s)$.\newline

When the action of $\s$ on $M$ is not free then the base space is
not anymore a smooth manifold, but a stratified pseudomanifold
$M/\s$ whose length depends on the number of orbit types. There is
a Gysin sequence of $M$ relating the DeRham cohomology of $M$ with
the intersection cohomology of $M/\s$
\[
    \cdots\rightarrow H^{^i}(M)\overset{\oint}{\rightarrow}
    H^{^{i-1}}_{_{\q-\overline{2}}}(M/\s)\overset{\varepsilon}\rightarrow
    H^{^{i+1}}_{_{\q}}(M/\s)\overset{\pi^{*}}\rightarrow
    H^{^{i+1}}(M)\rightarrow\cdots
\]
where $\q$, $\overline{2}$ are perversities in $M/\s$. The
connecting homomorphism is again the multiplication by the Euler
class $\varepsilon\in H^{^2}_{_{\overline{2}}}(M/\s)$. The fixed points' subspace
$M^{\s}$ is naturally contained in $M/\s$.
The link of a fixed stratum $S\subset M/\s$ is always
a cohomological complex projective space
\cite{S1},\cite{coloquio santiago}. \newline

In this article we extend the above situation for any
pseudomanifold $X$ and any action of $\s$
on $X$ preserving the local structure. The orbit map
$\pi:X\rightarrow X/\s$ induces a long exact sequence
\[
    \cdots\rightarrow
    H^{^i}_{_{\q}}(X)\rightarrow
    H^{^{i}}(\G(X/\s))\overset{\partial}\rightarrow
    H^{^{i+1}}_{_{\q}}(X/\s)\overset{\pi ^{*}}{\rightarrow}
    H_{_{\q}}^{^{i+1}}(X)
    \rightarrow\cdots
\]
relating the intersection cohomologies of $X$ and $X/\s$ with a
third term $H^{^{*}}(\G(X/\s))$ whose cohomology can be given in
terms of local and global basic cohomological data; we call it the
{\it Gysin term}. The above long exact sequence is the {\it Gysin
sequence}.\newline

Global data concerns the Euler class $\varepsilon\in H^{^2}_{_{\overline{2}}}(X/\s)$.
For instance, if $\varepsilon=0$ then $H^{^{*}}(\G(X/\s))=
H^{^{*}}_{_{\q-\overline{\chi}}}(X/\s)$
where $\overline{\chi}$ is the perversity defined by
\[
    \overline{\chi}(S)=\left\{
    \begin{array}{ll}
        1 & \text{$S$ a fixed stratum} \\
        0 & \text{else}
    \end{array}
    \right.
\]
The connecting homomorphism $\partial$ of the Gysin sequence depends on
the Euler class, though it's not the multiplication. The Euler class
vanishes if and only if there is a foliation on $X-\Sigma$
transverse to the orbits of the action \cite{JI1}, \cite{euler class}.\newline

Local data relates the Gysin term with the fixed strata. In
general, there is a second long exact sequence
\[
    \cdots
    \rightarrow
    H_{_{\q-\overline{\chi}}}^{^{i}}(X/\s)
    \rightarrow H^{^i}(\Resizq(X/\s))
    \overset{\partial'}\rightarrow
    H^{^{i+1}}(\G(X/\s))
    \overset{\imath^{*}}\rightarrow
    H_{_{\q-\overline{\chi}}}^{^{i+1}}(X/\s)
    \rightarrow\dots
\]
the residual term satisfying
\[
    H^{^{*}}(\Resizq(X/\s))=\underset{S}\prod\ H^{^{*}}
    (S,\mathfrak{Im}(\varepsilon_L))
\]
where $S$ runs over the fixed strata and
$H^{^{*}}(S,\mathfrak{Im}(\varepsilon))$ is the cohomology of $S$
with values on a locally trivial constructible presheaf
\cite{borel2} $\mathfrak{Im}(\varepsilon_L)$ with stalk
\[
    \mathcal{F}=\Im\{\varepsilon_L:H^{^{\q(S)-1}}(L/\s)
    \rightarrow H^{^{\q(S)+1}}(L/\s)\}
\]
the image of the multiplication by the Euler class
$\varepsilon_L\in H^{^2}(L/\s)$ of the action on the Link $L$ of
$S$. Since $L$ may not be a sphere, this term could not
vanish.\newline

\section*{Acknowledgments}
I would like to thank some helpful conversations with M. Saralegi
and F. Dalmagro, so as the accurate comments of the journal's
referee. While writing this article, I received the financial
support of the CDCH-Universidad Central de Venezuela and the Math
Department-Euskal Herriko Unibertsitatea, so as and the
hospitality of the staff in the Math Department-Universit\'e
D'Artois.


\begin{thebibliography}{99}


\bibitem{borel} BOREL, A. \& SPALTENSTEIN, N. "Sheaf theoretic intersection cohomology
(Bern, 1983)" in {\sl Intersection Cohomology-Swiss seminars.}
(Bern, 1983). Progress in Mathematics Vol.{\bf 50}, 47-182.
Birkh\"auser. Boston (1984).

\bibitem{bredon1}  BREDON, G. {\sl Introduction to Compact Transformation Groups.\/}
Pure and Applied Mathematics Vol.{\bf 46}. Academic Press. New
York (1972).

\bibitem{bott}  BOTT, R. \& LU T. {\sl Differential Forms in Algebraic Topology.\/}
Graduate Texts in Mathematics Vol.{\bf 82}. Springer-Verlag. New
York-Heidelberg- Berlin (1982).

\bibitem{brylinsky} BRYLINSKI, J.L. {\sl Equivariant Intersection
Cohomology}. {\tt Contemporary Math.} {\bf 132} (1992), 5-32.

\bibitem{dalmagro} DALMAGRO, F. {\sl Cohomolog\'{\i}a de intersecci\'on
de las acciones t\'oricas iteradas.} Ph.D. Thesis. Universidad
Central de Venezuela. Caracas (2003).

\bibitem{gm1}  GORESKY, M. \& MACPHERSON, R. {\sl Intersection Homology Theory.\/}
{\tt Topology} {\bf 19}, 135-162  (1980).

\bibitem{gm2} GORESKY, M. \& MACPHERSON, R. {\sl Intersection Homology
II.\/} {\tt Invent. Math.} {\bf 71}, 77-129  (1983).


\bibitem{borel2} HAEFLIGER, A. "Introduction to  piecewise linear
 intersection homology" (Bern, 1983)"
in {\sl Intersection Cohomology-Swiss seminars.} (Bern, 1983).
Progress in Mathematics Vol.{\bf 50}, 1-22. Birkh\"auser. Boston
(1984).

\bibitem{ghv}  GREUB, W.; HALPERIN, S. \& VANSTONE, R. {\sl Conections, curvature and
cohomology.\/} Pure and Applied Mathematics Vol.{\bf 47}.
Academic Press. New York (1972).

\bibitem{S1}  HECTOR, G. \& SARALEGI, M. {\sl Intersection Cohomology of
$\s$-actions.\/} {\tt Trans. Amer. Math. Soc.} {\bf 338}, 263-288
(1983).

\bibitem{king2} KING, H. {\sl Intersection homology
and homology of manifolds.\/} {\tt Topology} {\bf 21} (1982),
229-234.

\bibitem{king} KING, H. {\sl Topology invariance of intersection
homology without sheaves.\/} Topology Appl. {\bf 20} (1985),
149-160.

\bibitem{pervsheaves} MACPHEARSON, R.{\sl Intersection Homology and Perverse
Sheaves.\/} {Colloquium Lectures, Annual Meeting of the Amer.
Math. Soc.} San Francisco (1991).

\bibitem{coloquio santiago} MASA, X.; MACIAS E. \& ALVAREZ J.
{\sl Analysis and Geometry in Foliated Manifolds.\/} World
Scientific. Santiago de Compostela (1994).

\bibitem{nagase} NAGASE, M. {\sl $\L^2$-cohomology and intersection cohomology of stratified
spaces.\/} {\tt Duke  Math. J.} {\bf 50}, 329-368  (1983).

\bibitem{normalizer} PADILLA, G. {\sl On normal stratified pseudomanifolds.\/} (2002)
To appear in {\tt Extracta Math.}

\bibitem{pflaum} PFLAUM, M. {\sl Analytic and Geometric study of Stratified
Spaces.} Lecture Notes in Mathematics Vol.{\bf 1768}. Springer.
Berlin (2001).

\bibitem{JI1} ROYO, J. {\sl The Euler Class for Riemannian Flows},
{\tt C. R. Acad. Sci. Paris}, t.332, Serie I, pp. 45--50  (2001).

\bibitem{JI2} ROYO, J.{\sl The Gysin Sequence for Riemannian Flows},
{\tt Contemporary Mathematics} v.288, pp. 415--419  (2001).

\bibitem{toricas} SARALEGI, M. {\sl Cohoomologie d'Intersection des
Actions Toriques Simples.\/} {\tt Indag. Math.} {\bf 33}, 389-417
(1996).

\bibitem{S3} SARALEGI, M. {\sl A Gysin Sequence for Semifree Actions of $S^3$.\/}
{\tt Proc. Amer. Math. Soc.} {\bf 118}, 1335-1345 (1993).

\bibitem{illinois}  SARALEGI, M. {\sl Homological Properties of Stratified Spaces.\/}
{\tt Illinois J. Math.} {\bf 38}, 47-70 (1994).

\bibitem{euler class}  SARALEGI, M. {\sl The Euler class for flows of isometries.\/}
{\tt Research Notes in Math.} {\bf 131}, 25-28 (1989).

\bibitem{thom} THOM, R. {\sl Ensembles et morphismes
stratifi\'es\/}. Bull. Amer. Math. Soc. {\bf 75},240-284 (1969).





\end{thebibliography}
\end{document}